\documentclass[11pt]{article}

\usepackage{fullwidth}
\usepackage{color}
\usepackage[color,matrix,arrow]{xy}
\usepackage[top=1.5in, bottom=1.5in, left=1in, right=1in]{geometry}
\usepackage{amsgen}
\usepackage{amsmath}
\usepackage{amstext}
\usepackage{amsbsy}
\usepackage{amsopn}
\usepackage{amsfonts}
\usepackage{amssymb}
\usepackage{eepic}
\usepackage{graphicx}
\usepackage{epsf}
\usepackage{pstricks}
\xyoption{all}

\def\Box{\square}
\def\edge{\relbar\joinrel\relbar}

\def\mapright#1{\smash{\mathop{\longrightarrow}\limits^{#1}}}
\def\arr{\mapright{}}
\def\tra#1{\smash{\mathop{\mid\kern
-1pt\joinrel\relbar\joinrel\relbar}\limits^{*}_{#1}}}
\def\longtra#1{\smash{\mathop{\mid\kern
-1pt\joinrel\relbar\joinrel\relbar\joinrel\relbar}\limits^{*}_{#1}}}
\def\vlongtra#1{\smash{\mathop{\mid\kern
-1pt\joinrel\relbar\joinrel\relbar\joinrel\relbar\joinrel\relbar}\limits^{*}_{#1}}}
\def\vvlongtra#1{\smash{\mathop{\mid\kern
-1pt\joinrel\relbar\joinrel\relbar\joinrel\relbar\joinrel\relbar\joinrel\relbar}\limits^{*}_{#1}}}
\def\vvvlongtra#1{\smash{\mathop{\mid\kern
-1pt\joinrel\relbar\joinrel\relbar\joinrel\relbar\joinrel\relbar\joinrel\relbar\joinrel\relbar}\limits^{*}_{#1}}}
\def\etra#1{\smash{\mathop{\mid\kern
-1pt\joinrel\relbar\joinrel\relbar}\limits_{#1}}}

\def\Rw{\Rightarrow}
\def\oo{\overline}

\def\I{{\cal{I}}}

\def\N{\mathbb{N}}

\def\sji{\mbox{sji}}

\def\dim{\mbox{dim}}

\def\mindeg{\mbox{mindeg}}
\def\lat{\mbox{Lat}}
\def\LR{\mbox{LR}}
\def\cl{\mbox{Cl}}

\def\sub{\mbox{Sub}}

\def\flatx{\mbox{Fl}}

\def\H{{\cal{H}}}
\def\BoR{{\cal{BR}}}

\def\p{\varphi}

\def\inv{^{-1}}

\def\bi{\begin{itemize}}
\def\ei{\end{itemize}}
\def\beq{\begin{equation}}
\def\eeq{\end{equation}}

\def\AB#1{B({#1})}

\newtheorem{T}{Theorem}[section]
\newcommand{\bt}{\begin{T}}
\newcommand{\et}{\end{T}}
\newcommand{\ftd}{$\square$\end{T}}

\newtheorem{Proposition}[T]{Proposition}
\newcommand{\bp}{\begin{Proposition}}
\newcommand{\ep}{\end{Proposition}}
\newcommand{\fpd}{$\square$\end{Proposition}}

\newtheorem{Lemma}[T]{Lemma}
\newcommand{\bl}{\begin{Lemma}}
\newcommand{\el}{\end{Lemma}}
\newcommand{\fld}{$\square$\end{Lemma}}

\newtheorem{Corol}[T]{Corollary}
\newcommand{\bc}{\begin{Corol}}
\newcommand{\ec}{\end{Corol}}
\newcommand{\fcd}{$\square$\end{Corol}}

\newtheorem{Result}[T]{Result}
\newcommand{\br}{\begin{Result}}
\newcommand{\er}{\end{Result}}
\newcommand{\frd}{$\square$\end{Result}}

\newtheorem{Example}[T]{Example}
\newcommand{\be}{\begin{Example}}
\newcommand{\ee}{\end{Example}}

\newtheorem{Problem}[T]{Problem}
\newcommand{\bq}{\begin{Problem}}
\newcommand{\eq}{\end{Problem}}

\newtheorem{Remark}[T]{Remark}
\newcommand{\brm}{\begin{Remark}}
\newcommand{\erm}{\end{Remark}}

\newcommand{\proof}
   {\par\medbreak\noindent{\bf Proof}.\enspace}

\newcommand{\qed}{
$\Box$
\par\bigbreak}

\normalbaselines
\def\abstract#1{\par\bigskip
\begingroup\small
\baselineskip=12truept
\begin{center}ABSTRACT\end{center}
\par\medskip\par\noindent
\null\hfill\hbox{\vbox{\hsize=5truein\noindent#1}}
\hfill\null\par\endgroup\par}

\title{On the subsemigroup complex of an aperiodic Brandt semigroup}
\author{{\bf Stuart Margolis, John Rhodes and Pedro V. Silva}}

\date{\today}

\begin{document}
\maketitle

\begin{center}\small
2010 Mathematics Subject Classification: 20M18, 05B35, 05E45, 55U10

\bigskip

Keywords: Brandt semigroup, lattice of subsemigroups, simplicial complex, boolean representable simplicial complex, matroid
\end{center}

\abstract{We introduce the subsemigroup complex of a finite semigroup $S$ as a (boolean representable) simplicial complex defined through chains in the lattice of subsemigroups of $S$. We present a research program for such complexes, illustrated through the particular case of combinatorial Brandt semigroups. The results include alternative characterizations of faces and facets, asymptotical estimates on the number of facets, or establishing when the complex is pure or a matroid.
}

\section{Introduction}

In a recent paper \cite{CGMP}, Cameron, Gadouleau, Mitchell and
Peresse discuss the maximal length of chains of subsemigroups for
various classes of semigroups. This problem can be viewed as the {\em
dimension problem} for a boolean representable simplicial complex
naturally associated with the lattice of subsemigroups of a
semigroup. It follows that there exist a number of natural questions
associated with this complex, which may shed a new light on the nature
of $S$. The present paper lays the foundations for the subsemigroup complex of a semigroup $S$
and then illustrates this program by
considering the lattice and complex of subsemigroups of an aperiodic Brandt semigroup.
As it turns out, this particular case suffices to
expose quite a number of astonishing connections. For the non aperiodic case, the reader is referred to another paper by the same authors \cite{MRS2}. The reader is assumed to be familiar with the foundations of semigroup theory (see \cite{RS}).

In a series of three papers \cite{IR1,IR2,IR3}, Izhakian and Rhodes
introduced the concept of boolean representation for various algebraic
and combinatorial structures. These ideas were inspired by previous
work by Izhakian and Rowen on supertropical matrices (see
e.g. \cite{Izh2,IRow2,IRow,IRow22}), and were subsequently developed by
Rhodes and Silva in a recent monograph, devoted to boolean
representable simplicial complexes \cite{RSm}. Note that simplicial
complexes may be viewed under two
perspectives, geometric and combinatorial. It is well known that each
structure determines the other (see e.g. \cite[Section A.5]{RSm}).

The original approach to boolean representable simplicial complexes is to consider matrix representations over the
superboolean semiring $\mathbb{SB}$, using appropriate notions of vector
independence and rank. Writing $\mathbb{N} = \{ 0,1,2,\ldots \}$, we
can define $\mathbb{SB}$ as the quotient of $(\mathbb{N},+,\cdot)$ (usual
operations) by the congruence which identifies all integers $\geq 2$.
In this context, boolean representation refers to matrices using only
0 and 1 as entries.

Equivalently, boolean representable simplicial complexes
can be characterized by
means of chains in lattices, namely in the lattice of flats. The
lattice of flats plays a fundamental
role in matroid theory but is not usually considered for arbitrary
simplicial complexes, probably due to the fact that, unlike in the
matroid case, the structure of a simplicial complex cannot be in general
recovered from its lattice of flats. However, this is precisely what
happens with boolean representable simplicial complexes. If $\H\; =
(A,H)$ is a simplicial complex and $\flatx\H$ denotes its lattice of
flats, then $\H$ is boolean representable if and only if $H$ equals
the set of transversals of the successive differences for chains in
$\flatx\H$ \cite{RSm}. This implies in particular that all (finite) matroids are boolean
representable. And this generalizes Birkhoff's Theorem \cite{Oxl} that gives a 1-1 correspondence between geometric lattices and simple matroids, to arbitrary lattices and simple boolean representable simplicial complexes.

\section{Boolean representable simplicial complexes}

All lattices, simplicial complexes and semigroups in this paper are
assumed to be
finite. Given a set $A$ and $n \geq 0$, we denote by $P_n(A)$
(respectively $P_{\leq n}(A)$) the set
of all subsets of $A$ with precisely (respectively at most) $n$
elements.

A (finite) simplicial complex is a structure of the form $\H\; =
(A,H)$, where $A$ is a finite nonempty set and $H \subseteq 2^A$
contains $P_1(A)$ and is closed under taking subsets. The elements of
$A$ and $H$ are called respectively {\em vertices} and {\em faces}.

A face of $\H$ which is maximal with respect to
inclusion is called a {\em facet}. We denote by ${\rm fct}\H$ the set
of facets of $\H$.

The {\em dimension} of a face $I \in H$ is $|I|-1$. An $i$-{\em face}
(respectively $i$-{\em facet}) is a face (respectively facet) of
dimension $i$. We may refer to 0-faces and 1-faces as vertices and edges.


We say that $\H$ is:
\bi
\item {\em simple} if $P_2(A) \subseteq H$;
\item {\em pure} if all the facets
of $\H$ have the same dimension.
\ei
The dimension of $\H$, denoted by $\dim\H$, is the maximum dimension of a
face (or facet) of $\H$.


Two simplicial complexes $(A,H)$ and $(A',H')$ are {\em isomorphic} if
there exists a bijection $\p:A \to A'$ such that
$$X \in H \mbox{ if and only if }X\p \in H'$$
holds for every $X \subseteq A$.

If $\H\; = (A,H)$ is a simplicial complex and $W \subseteq A$ is
nonempty, we call
$$\H|_W = (W,H \cap 2^W)$$
the {\em restriction} of $\H$ to $W$. It is obvious that $\H|_W$ is
still a simplicial complex.

A simplicial complex $\H\; = (A,H)$ is called a {\em matroid} if it
satisfies the {\em exchange property}:
\bi
\item[(EP)]
For all $I,J \in H$ with $|I| = |J|+1$, there exists some
  $p \in I\setminus J$ such that $J \cup \{ p \} \in H$.
\ei

A simplicial complex $\H\; = (A,H)$ is {\em
  shellable} if we can order its facets as $B_1, \ldots, B_t$ so
that, for $k = 2, \ldots,t$ and if
$I(B_k) = (\cup_{i=1}^{k-1} 2^{B_i}) \cap 2^{B_k}$, then
$$(B_k,I(B_k))\mbox{ is pure of dimension }|B_k|-2$$
whenever $|B_k| \geq 2$.
Such an ordering is
called a {\em shelling}.
Shellability is an important property since it implies that
the complex has the homotopy type of a wedge of spheres \cite{BW}.

Given an $R \times A$ matrix $M$ and $Y \subseteq R$, $X \subseteq A$,
we denote by $M[Y,X]$ the submatrix of $M$ obtained by deleting all
rows (respectively columns) of $M$ which are not in $Y$ (respectively $X$).

A boolean matrix $M$ is {\em lower unitriangular} if it is of the form
$$\left(
\begin{matrix}
1&&0&&0&&\ldots&&0\\
?&&1&&0&&\ldots&&0\\
?&&?&&1&&\ldots&&0\\
\vdots&&\vdots&&\vdots&&\ddots&&\vdots\\
?&&?&&?&&\ldots&&1
\end{matrix}
\right)
$$

Two matrices are {\em congruent} if we can transform one into the other
by independently permuting rows/columns. A boolean matrix is {\em nonsingular}
if it is congruent to a lower unitriangular matrix.

Given an $R \times A$ boolean matrix $M$, we say that the
subset of columns $X \subseteq A$ is $M$-{\em independent} if there exists
some $Y \subseteq R$ such that $M[Y,X]$ is nonsingular.

A simplicial complex $\H\; = (A,H)$ is {\em boolean representable} if
there exists some boolean matrix $M$ such that
$H$ is the set of all $M$-independent subsets of $A$. We denote by
$\mindeg(\H)$ the minimum number of rows of such a matrix (if it exists).

We shall use the acronym BRSC to denote a (finite) boolean representable simplicial
complex, and we denote by $\BoR$
the class of all BRSCs. All matroids are boolean representable, but the converse is
not true.

We say that $X
\subseteq A$ is a {\em flat} of $\H$ if
$$\forall I \in H \cap 2^X \hspace{.2cm} \forall p \in A \setminus X
\hspace{.5cm} I \cup \{ p \} \in H.$$
The set of all flats of $\H$ is denoted by
$\flatx\H$. Note that $A, \emptyset \in \flatx\H$ in all cases.

Clearly, the intersection of any set of flats (including $A =
\cap\emptyset$) is still a flat. If we order $\flatx\H$ by inclusion,
it is then a $\wedge$-semilattice. Since $\flatx\H$ is finite and contains a maximal element,
it is  a lattice for the determined join, that is, the join of two flats is the intersection of all flats containing their union.
We call $\flatx\H$ the {\em lattice of flats} of
$\H$. The lattice $\flatx\H$ induces a closure operator on $2^A$ defined by
$$\oo{X} = \cap\{ F \in \flatx \H \mid X \subseteq F \}$$
for every $X \subseteq A$.

The lattice of flats is an important example of a lattice which is
$\vee$-generated by $A$. A lattice $L$ is $\vee$-{\em generated} by $A$
if there exists a mapping $\iota:A \to L$ such that
$$L = \{ \vee(B\iota) \mid B \subseteq A \}.$$
We assume that $\vee \emptyset$ is the bottom element of
$L$ by convention.
In the case of $\flatx\H$, the mapping $\iota:A \to \flatx\H$ is
predictably defined by $a\iota = \oo{a}$ $(a \in A)$.

We shall represent structures of the above type as
ordered pairs $(L,\iota)$. The class of all $(L,\iota)$ will be
denoted by $\lat_{\vee}(A)$.

Let $(L,\iota) \in \lat_{\vee}(A)$. We say that $X \subseteq A$ is a
{\em transversal of the successive differences} for a chain
$$\ell_0 < \ell_1 < \ldots < \ell_k$$
if $X$ admits an enumeration $x_1,\ldots , x_k$ such that $\ell_{i-1}
< (\ell_{i-1} \vee x_i\iota)
\leq \ell_i$ for $i = 1,\ldots,k$. We say that $(L,\iota)$ is a {\em
  lattice representation} of a simplicial complex $\H\; = (A,H)$ if the
elements of $H$ are precisely the transversals of the successive
differences for chains in $L$. Conversely, the set of all
transversals of the successive differences for chains in $(L,\iota)
\in \lat_{\vee}(A)$ constitutes a boolean representable simplicial
complex $\H(L,\iota) = (A,H(L,\iota))$.

We denote by $\LR(\H)$ the class of all
lattice representations of $\H$. A simplicial complex is boolean
representable if and only if it admits a lattice representation
\cite[Section 5.4]{RSm}. Up to isomorphism, every such lattice
representation may be viewed as a sublattice of $\flatx\H$, which
plays then the canonical role of being the largest lattice
representation.
Moreover, by \cite[Corollary 5.2.7]{RSm}, $\H\; = (A,H)$ is boolean
representable if and only if
\beq
\label{derby3}
\mbox{every $X \in H$ admits an enumeration
$x_1,\ldots, x_k$ satisfying $\oo{x_1} \subset \oo{x_1x_2}
  \subset\ldots \subset \oo{x_1\ldots x_k}.$}
\eeq

Now we may define an ordering on $\lat_{\vee}(A)$
by $(L,\iota) \geq (L',\iota')$ if there exists a $\vee$-map
(i.e. preserving arbitrary joins) $\p:L \to L'$ such that the diagram
$$\xymatrix{
&A \ar[dl]_{\iota} \ar[dr]^{\iota'} & \\
L \ar[rr]_{\p} && L'
}$$
commutes. This is an appropriate comma category and we quotient by the equivalence relation that identifies two lattice representations if they are related to one another by $\leq$. We thus obtain a partially ordered set, or equivalently a skeletal category with at most one morphism between two objects.
Then $\LR(\H)$ constitutes an up set of $\lat_{\vee}(A)$ and we may
identify all the elements of $\lat_{\vee}(A) \setminus \LR(\H)$ into a
single bottom element $B$ to obtain $\LR_0(\H) = \LR(\H) \cup \{ B
\}$. If $\rho = (\leq \cap \geq)$, the quotient $\LR_0(\H) /\rho$
constitutes a lattice, the {\em lattice of lattice representations} of
$\H$. Note that the bottom element is not a true lattice
representation (it agglomerates instead the non representations), but
the atoms of this lattice (the {\em minimal representations} of $\H$)
are the most economical ways of representing $\H$ through a
lattice. The strictly join irreducible (sji) elements of the lattice of
lattice representations (called join irreducible in many lattice books) are also important, see \cite[Chapter 5]{RSm}
for details. An element $a$ of a lattice $L$ is {\em sji} if, for
every $X \subseteq L$, $a = \vee X$ implies $a \in X$. This is
equivalent to saying that $a$ covers exactly one element of $L$. We
denote by $\sji(L)$ the set of all sji elements of $L$.




\section{Simplification}
\label{simply}

Let $\H\; = (A,H) \in \BoR$. Recall that $\H$ is simple if every pair of distinct elements is a face of $\H$. In this section we show that the simplification
of an $\H\; = (A,H) \in \BoR$ is also boolean representable. This generalizes a well known result in matroid theory.

We define an equivalence relation $\eta$
on $A$ by
$$a \eta b \; \mbox{ if } \; \oo{a} = \oo{b}.$$
It is easy to see that if $\H = \H(L,\iota)$ for some $(L,\iota) \in
\lat_{\vee}(A)$, then $a \eta b$ if and only if $a\iota =
b\iota$. Indeed, both conditions are clearly equivalent to $\{
a,b\} \notin H$.

Now we define a
simplicial complex $\H_{\eta} = (A/\eta,H/\eta)$, where
$$H/\eta = \{
\{ a_1\eta, \ldots,a_k\eta \} \mid \{ a_1,\ldots,a_k \} \in H \}.$$
Note that, in view of (\ref{derby3}),
\beq
\label{retro1}
\mbox{if the $a_i$ are all distinct, so must be the
$a_i\eta$.}
\eeq
\beq
\label{retro2}
\mbox{if $a_i\eta b_i$ for $i = 1,\ldots,k$, then
$\{ a_1,\ldots,a_k \} \in H$ if and only if $\{ b_1,\ldots,b_k \} \in
H$.}
\eeq
It follows that $\H_{\eta}$ is isomorphic to the restriction
$\H|_W$ for any cross-section $W$ of $\eta$. We call $\H_{\eta}$ the {\em simplification} of $\H$.

We collect in the next result some of the properties of the simplification. All claims were proved in \cite[Proposition 4.2]{MRS} except for (viii) and (ix), which we prove below.

\bp
\label{retro}
Let $\H\; = (A,H) \in \BoR$ and let $\p:A \to
A/\eta$ be the canonical projection. Then:
\bi
\item[(i)] ${\rm dim}\H/\eta = {\rm dim}\H$;
\item[(ii)] ${\rm Fl}\H\; = \{ F\p\inv \mid F \in {\rm Fl}(\H/\eta)
  \}$;
\item[(iii)] ${\rm Fl}\H\cong {\rm Fl}(\H/\eta)$;
\item[(iv)] $\H/\eta$ is boolean representable;
\item[(v)] $\H/\eta$ is simple;
\item[(vi)] $\H$ is pure if and only if $\H/\eta$ is pure;
\item[(vii)] $\H$ is a matroid if and only if $\H/\eta$ is a matroid;
\item[(viii)] ${\rm mindeg}(\H) = {\rm mindeg}(\H/\eta)$;
\item[(ix)] the lattices of lattice representations of $\H$ and
  $\H/\eta$ are isomorphic;
\item[(x)] if $\H/\eta$ is shellable, so is $\H$.
\ei
\ep

\proof
(viii) Assume that $M$ is an $R \times A$ boolean matrix representing
$\H$. Since $P_1(A) \subseteq H$, $M$ has no zero columns. On the
other hand, if $\oo{a} = \oo{b}$, then $\{ a,b \} \notin H$ by
(\ref{derby3}) and so the $a$th and $b$th columns of $M$ are equal.

It follows easily from (\ref{retro1}) and (\ref{retro2}) that every
matrix representation of $\H$ induces a matrix representation of
$\H/\eta$ by removing repeated columns inside each
$\eta$-class. Conversely, every
matrix representation of $\H\eta$ induces a matrix representation of
$\H$ by adding repeated columns for the redundant elements inside each
$\eta$-class. Therefore the minimum possible degree must be the same.

(ix) The lattice representations of $\H$ correspond to all quotients
of $(\flatx\H,\iota)$ which retain the capacity of recognizing all the
elements of $H$ as transversals of the successive differences for
chains.

Given $X \subseteq A$, let $\cl_{\eta}(X\eta)$ denote the closure of $X\eta$ in $\H_{\eta}$.

By part (ii), the lattices $\flatx\H$ and $\flatx(\H/\eta)$
are isomorphic, and the isomorphism is compatible with the image of
the generating sets, $a \mapsto \oo{a}$ and $a\p \mapsto
\cl_{\eta}(a/\eta)$ $(a \in A)$. On the other hand, in view of
(\ref{retro1}) and (\ref{retro2}), we must test essentially the same
chains. We omit the technical details, but we are indeed led to a
canonical isomorphism between the lattices of lattice representations
of $\H$ and $\H/\eta$.
\qed

\section{The subsemigroup complex of a semigroup}

Let $S$ denote a finite nonempty semigroup and let $\sub(S)$ denote the lattice
of subsemigroups of $S$, ordered under inclusion. To have a lattice,
we view the empty set as a semigroup. Let $\iota_S:S \to \sub(S)$ be
defined by $s\iota_S = s^+$, the subsemigroup generated by $\{s\}$. Then $(\sub(S),\iota_S) \in
\lat_{\vee}(S)$. Write $H(S) = H(\sub(S),\iota_S)$.
The {\em subsemigroup complex} of $S$ is the boolean
representable simplicial complex $\H(S) = (S,H(S)).$
That is, a subset $X \subseteq S$ is a face if and only if it admits
an enumeration $x_1,\ldots,x_n$ such that, for some chain
$$S_0 \subset S_1 \subset \ldots \subset S_n$$
of subsemigroups of $S$, we have $x_i \in S_i \setminus S_{i-1}$ for $i =
1,\ldots,n$. If such a chain exists, it can be taken as
\beq
\label{hsc}
\emptyset \subset x_1^+ \subset \{ x_1,x_2\}^+ \subset \ldots
\subset \{ x_1,\ldots, x_n\}^+.
\eeq

In general, $\H(S)$ is not simple since different elements can
generate the same subsemigroup.

\bl
\label{etahs}
Let $S$ be a finite nonempty semigroup and consider $\H(S) \in
\BoR$. For all $s,t \in S$,
$$\mbox{$\oo{s} = \oo{t}$ if and only if $s^+ = t^+$}.$$
\el

\proof
By (\ref{derby3}), we have $\oo{s} = \oo{t}$ if and only if $\{ s,t \}
\notin H(S)$. This is equivalent to avoiding chains of the form
(\ref{hsc}), i.e. $s^+ = t^+$.
\qed

We define also $\H_0(S) = (S/\eta,H_0(S)) = \H(S)/\eta$. We may
identify $S/\eta$ with the set of cyclic subsemigroups of $S$. A set
$Y$ of cyclic subsemigroups is a face if and only if it admits an
enumeration $C_1,\ldots,C_n$ such that
$$C_1 \subset (C_1 \cup C_2)^+ \subset \ldots \subset (C_1 \cup \ldots
\cup C_n)^+.$$
In view of Proposition \ref{retro}, the complexes $\H(S)$ and
$\H_0(S)$ are equivalent with respect to most properties, we can work
with either at our convenience.

Now we note a property about finite
semigroups all of whose nontrivial subgroups have prime order. This includes
the case of finite {\em aperiodic} semigroups (all subgroups are
trivial, or equivalently, satisfying an identity
$x^{n+1} = x^n$ for some $n \in \N$).

A finite lattice $L$ has a minimum
$\vee$-generating set, namely $\sji(L)$. Let $\mu_L:\sji(L)
\to L$ denote the inclusion map. Then $(L,\mu_L) \in
\lat_{\vee}(\sji(L))$ and so $\H(L,\mu_L) \in \BoR$.

\bp
\label{ape}
Let $S$ be a finite semigroup whose nontrivial subgroups have prime
order. Then $\H_0(S) \cong \H({\rm Sub}(S),\mu_{{\rm
    Sub}(S)})$. Moreover, if $S$ is aperiodic, then $\H(S) = \H_0(S)$.
\ep

\proof
For the first claim, it suffices to show that the cyclic subsemigroups
of $S$ are precisely
the sji elements of $\sub(S)$. Since any subsemigroup of $S$ is a join of
cyclic subsemigroups, every sji subsemigroup must be necessarily
cyclic. Conversely, let $s \in S$ and suppose that $s^+ = (s_1^+ \vee
\ldots \vee s_k^+) = \{ s_1,\ldots, s_k\}^+$. Suppose first that $s^+$
is not a subgroup. Then $ss^+ = s^+ \setminus \{ s \}$ is a
subsemigroup of $S$. If $s \notin \{
s_1,\ldots, s_k\}$, then $\{ s_1,\ldots, s_k\}^+
\subseteq ss^+ \in \sub(S)$, a contradiction. Hence
$s \in \{ s_1,\ldots, s_k\}$ and so $s^+$ is an sji as required.

Suppose now that $s^+$ is a subgroup. If the order of $s$ (i.e. $|s^+|$) is $1$, then $s_1 =
\ldots s_k = s$ and we are done. Hence we may assume that $|s^+|$ is a
prime. Then $s^+ = \{ s_1,\ldots, s_k\}^+$ implies that $s_i$ is not
the identity for some $i$, thus $s_i^+ = s^+$ and so $s^+$ is an sji
as required.

If $S$ is aperiodic, then $s$ is the unique generator of $s^+$, whence
$\eta$ is the identity and so $\H_0(S) = \H(S)$.
\qed

As a consequence, if $S$ is a finite semigroup whose nontrivial
subgroups have prime order,
we can say that the properties
of $\H(S)$ are fully determined by the lattice $\sub(S)$. Note that
Proposition \ref{ape} does not hold for arbitrary finite semigroups:
the cyclic group $C_6$ provides an easy counterexample.

We can now enumerate a list of very natural problems which can be
stated in connection with the boolean representable simplicial
complexes $\H(S)$ and $\H_0(S)$:
\bi
\item[(P1)]
To compute dimension.
\item[(P2)]
To characterize the faces.
\item[(P3)]
To characterize the flats.
\item[(P4)]
To compute the lattice of flats.
\item[(P5)]
To determine whether the complex is pure.
\item[(P6)]
To determine whether the complex is a matroid.
\item[(P7)]
To determine whether the complex is shellable.
\item[(P8)]
To compute the minimum degree of a matrix representation.
\item[(P9)]
To compute the minimal lattice representations.
\ei

\section{Subsemigroups of aperiodic Brandt semigroups}

\subsection{Basic properties of $B(n)$ and its subsemigroups}

In this section we study subsemigroups of the aperiodic Brandt semigroup $\AB{n}$. The semigroup $\AB{n}$ can be defined as the set $(\{1,...,n\} \times \{1,...,n\}) \cup \{0\}$ where 0 is the zero element and product $(i,j)(k,l) = \delta(j,k)(i,l)$, where $\delta(.,.)$ is the Kronecker delta. $\AB{n}$ arises in a number of different ways in mathematics that we look at in detail. The interaction between the various ways of thinking about $\AB{n}$ and its subsets gives us the opportunity to deepen our understanding of these objects. If $X$ is an arbitrary set, we will use the notation $\AB{X}$ for the corresponding aperiodic Brandt semigroup on the set $(X \times X) \cup \{0\}$.

First of all, $\AB{n}$ is a 0-simple inverse semigroup. The inverse of $(i,j)$ is $(j,i)$. As a Rees matrix semigroup, \cite{CP}, $\AB{n}$ is isomorphic to the Rees matrix semigroup $M^{0}(\{1,...,n\}, 1, \{1,...,n\},I_{n})$, where $I_{n}$ is the $ n \times n$ identity matrix. $\AB{n}$ and its inverse subsemigroups (detailed below) play an important role in inverse semigroup theory.

Secondly, if we restrict the multiplication of $\AB{n}$ to all products that are non-zero, then $\AB{n}$ is a groupoid (in the sense of category theory, a category all of whose morphisms are isomorphisms.) It is the unique connected trivial groupoid with $n$ objects. This means that there is exactly one morphism between any two objects. We'll see below that we can identify the subsemigroupoids (similar to categories, but may not have an identity at each object) of $\AB{n}$ where $n$ runs over the natural numbers, with the collection of all finite trivial semigroupoids.

Thirdly, if we identify the pair $(i,j)$ of $\AB{n}$ with the elementary $n \times n$ matrix which has entry 1 in position $(i,j)$ and 0 in all other positions and 0 with the 0 matrix, then we can identify $\AB{n}$ with the semigroup of elementary $n \times n$ matrices over any non-trivial semiring. Since we are interested in subsets of $\AB{n}$ it is natural to look at the monoid $M_{n}(B)$ of all $n \times n$ matrices over the two element Boolean semiring.

Finally, if we think of the element $(i,j)$ of $\AB{n}$ as an ordered pair, we see that non-zero elements of $\AB{n}$ can be thought of as binary relations with one element, or equivalently, as directed graphs with one edge. Subsets of $\AB{n}$ and in particular subsemigroups of $\AB{n}$ can then be uniquely identified with arbitrary binary  relations on $\{1,...,n\}$ or equivalently as directed graphs with vertex set $\{1,...,n\}$. This allows us to use tools from graph theory to study subsemigroups of $\AB{n}$.

Notice that the semigroup $P(\AB{n})$ of all  subsets of $\AB{n}$ is a monoid with identity element $1=\{(i,i)|i=1, \ldots, n\}$.The next theorem is a straightforward calculation whose proof is left to the reader.

\bt
\label{Rels} Let $P(\AB{n})$  be the monoid of all subsets of $\AB{n}$ and $R_{n}$ the monoid of all binary relations on $\{1, \ldots, n\}$. The function f:$P(\AB{n}) \rightarrow R_n$ defined by
$f(X) =  X \setminus \{ 0 \}$
is a surjective morphism.
\et


The morphism $f$ in Theorem \ref{Rels} is exactly 2 to 1, only identifying $X$ with $X \cup \{0\}$. The subset $P_{0}(\AB{n})$ of $P(\AB{n})$ consisting of all subsets containing 0 is a subsemigroup of $P(\AB{n})$ and a monoid with identity $1=\{(i,i)|i=1 \ldots n\} \cup \{0\}$. Thus the restriction of the morphism $f$ in Theorem \ref{Rels} is an isomorphism from $P_{0}(\AB{n})$ to $R_{n}$. Furthermore, it is well known that $R_{n}$ is isomorphic to the monoid $M_{n}(B)$ of all $n \times n$ Boolean matrices.

We now turn to subsemigroups of $\AB{n}$. Clearly the only subsemigroups of $\AB{n}$ that do not contain 0, are the empty semigroup, and all one element semigroups $\{(i,i)\}$ for $ i \in \{1,\ldots, n\}$. In the following discussion, when we speak of subsemigroups of $\AB{n}$, we mean those that contain 0. We wish to exploit the isomorphism in Theorem \ref{Rels} to get an interpretation of subsemigroups within the context of binary relations and directed graphs.

Since subsemigroups of a semigroup $S$ are precisely the subsets $T$ of $S$ such that $T^{2} \subseteq T$, it follows that under the isomorphism from $P_{0}(\AB{n})$ to $R_{n}$, the subsemigroups of $\AB{n}$ containing 0 can be identified with the transitive relations on $\{1, \ldots, n\}$, which are by definition the relations $R$ such that $R^{2} \subseteq R$.

Inverse subsemigroups of $\AB{n}$ are the subsemigroups that are closed under the operation that sends $(i,j)$ to $(j,i)$. Under the identification as binary relations above, inverse subsemigroups are exactly the transitive and symmetric relations. These in turn are precisely the partial equivalence relations on $\{1, \ldots, n\}$, that is, an equivalence relation on a subset of $\{1, \ldots, n\}$. Thus an inverse subsemigroup $I$ of $\AB{n}$ can be identified with a partition of the set $\{i|(i,i) \in I\}$.

A subsemigroup $T$ of a semigroup $S$ is called a full subsemigroup if $E(T) = E(S)$, where $E(S)$ is the set of idempotents of $S$. As relations it is then clear that full subsemigroups of $\AB{n}$ correspond to reflexive and transitive relations, that is, they are precisely the preorders on $\{1, \ldots, n\}$. It is well known that this is the same as the set of topologies on an $n$ element set. The preorder associated to a topology $\mathcal{T}$ is the relation $x R y$ if and only if $x$ belongs to every neighborhood of $y$. If we think of $\AB{n}$ as the trivial connected groupoid on $n$ objects, then the preorders are exactly its subcategories, which are precisely the trivial categories on $n$ objects.

Finally, full inverse subsemigroups of $\AB{n}$ correspond to equivalence relations, a fact first noted by Jones \cite{Jones1}. We record the above observations in the
following theorem. We let \\
($F\sub_{0}(\AB{n}), I\sub_{0}(\AB{n}), FI\sub_{0}(\AB{n})$) $\sub_{0}(\AB{n})$ denote the lattice of (full, inverse, full inverse) subsemigroups of $\AB{n}$ containing 0.

\bt \label{subs}

\begin{itemize}

 \item[(i)] $\sub_{0}(\AB{n})$ is isomorphic to the lattice of transitive relations on an $n$ element set.

  \item[(ii)] $F\sub_{0}(\AB{n})$ is isomorphic to the lattice of preorders on an $n$ element set. Equivalently, it is isomorphic to the lattice of topologies on an $n$ element set.

  \item[(iii)] $I\sub_{0}(\AB{n})$ is isomorphic to the lattice of partial partitions on an $n$ element set.

  \item[(iv)] $FI\sub_{0}(\AB{n})$ is isomorphic to the lattice of partitions on an $n$ element set.

\end{itemize}

\et

It is well known that the lattice of partial partitions on an $n$ element set is isomorphic to the lattice of partitions on a set of size $n+1$. More formally, let $\Pi$  be a partition on $\{1, \ldots n+1\}$, and let $B$ be the block of $\Pi$ that contains $n+1$. Then $\Pi \backslash \{B\}$ is a partial partition on $\{1, \ldots n\}$. Conversely, if $\Phi$ is a partition on a subset $Y$ of $\{1, \ldots n\}$, we can define the partition on $\{1, \ldots n+1\}$ whose blocks are those $\Phi$ and the one extra block $(\{1, \ldots n\}\backslash Y) \cup \{n+1\}$. These operations are easily seen to be lattice isomorphisms inverse to one another. Therefore, $I\sub_{0}(\AB{n})$ and $FI\sub_{0}(\AB{n})$ are geometric lattices as this is true of partition lattices. On the other hand, it is easily seen that the lattice of transitive relations on a set of size at least 2 is not a semimodular lattice and thus is not a geometric lattice.

Theorem \ref{subs} allows us to count the various types of semigroups mentioned there by using results on the corresponding type of relation. For example, the number of subsemigroups of $\AB{n}$ is equal to $T(n) + n + 1$, where $T(n)$ is the number of transitive relations on an $n$ element set. The extra $n+1$ accounts for those subsemigroups not containing 0 (including the empty subsemigroup). There has been a good deal of literature on problems related to counting these objects. See \cite{Klaska, Pfeiffer} for example.

We now look at the structure of the various subsemigroups of $\AB{n}$ via the correspondences given in Theorem \ref{subs}. If $S$ is a subsemigroup of $\AB{n}$, we let $\Gamma(S)$ denote the directed graph corresponding to the (transitive) relation corresponding to $S$ given by Theorem \ref{subs}.

We first look at inverse subsemigroups. Recall that if $\{S_{i}| i \in I\}$ is a collection of semigroups, possibly with 0, then their 0-disjoint union is the disjoint union
$S = \{S_{i} \backslash \{0\}| i \in I\} \cup \{0\}$ with the product given by that of $S_i$ within each component (and when the product is 0 in $S_i$, it becomes 0 in $S$) and the product of elements in different components equal
to 0. The following theorem is straightforward to prove and we leave its proof to the reader.

\bt \label{invsub}

Let $\Pi = \{X_{i}| i \in I\}$ be a partial partition of $\{1, \ldots , n\}$. If $|I| > 1$, then the inverse subsemigroup of $\AB{n}$ corresponding to $\Pi$ is
the 0-disjoint union of $\{\AB{X_{i}}| i \in I\}$.

\et

We now look at the case of nilpotent subsemigroups. Recall that a nilpotent semigroup is either the empty semigroup or a semigroup $N$ with 0, such that $S^{k}= 0$ for some positive integer $k$. That is, the product of any $k$ elements of $N$ is 0. The least such integer $k$ is called the index of nilpotency of $N$. We assume below that we are talking about nonempty nilpotent semigroups. Finite nonempty nilpotent semigroups are easily seen to be the finite semigroups that have a unique idempotent that is the 0 element.

\bt \label{nilpsub}

Let $N$ be a subsemigroup of $\AB{n}$ containing 0. Then the following are equivalent.

\begin{itemize}

  \item[(i)] $N$ is a nilpotent semigroup.

  \item[(ii)] As a binary relation, $N$ is an irreflexive transitive relation. That is, $N$ is a strict partial order.

  \item[(iii)] The graph $\Gamma(N)$ is an acyclic directed graph.

\end{itemize}
\et

\proof Let $N$ be a nilpotent subsemigroup of $\AB{n}$ containing 0. Then $(i,i) \notin N$ for all $1 \leq i \leq n$ since 0 is the unique idempotent of $N$. $N$ is transitive since it is a subsemigroup. Therefore (i) implies (ii)

If $\Gamma(N)$ contains a cycle, from $i$ to $i$, then the idempotent $(i,i) \in N$ and $N$ is not an irreflexive relation and (ii) implies (iii). Finally, if $\Gamma(N)$ is a directed acyclic graph, then there is a largest integer $k$ such that $\Gamma(N)$ has a path of length $k$. It follows that the product of any $k+1$ elements of $N$ is 0 and $N$ is a nilpotent semigroup of nilpotency index $k+1$. \qed

We note that the intersection of nilpotent subsemigroups of $\AB{n}$ is a nilpotent subsemigroup, but the join in the lattice of all subsemigroups of $\AB{n}$ need not be nilpotent. For example the join of the two nilpotent subsemigroups $\{(1,2), 0\}$ and $\{(2,1), 0\}$ of $\AB{2}$ is all of $\AB{2}$. It is clear that maximal nilpotent subsemigroups of $\AB{n}$ correspond by Theorem \ref{nilpsub} to strict linear orders. Thus there are precisely $n!$ maximal nilpotent subsemigroups of $\AB{n}$, one for each way of listing all the elements of $\{1, \ldots, n\}$. This is a special case of a result of Graham \cite{Graham}, who classified the maximal nilpotent subsemigroups of an arbitrary finite 0-simple semigroup by graph theoretic methods. The problem of counting the number of maximal nilpotent subsemigroups containing a given nilpotent subsemigroup is thus the same as counting the number of linear extensions of a given partial order, a well studied and computationally difficult problem.

We can now describe the structure of arbitrary subsemigroups of $\AB{n}$. In the following theorem we consider the empty semigroup to be both a nilpotent semigroup and an inverse semigroup.

\bt \label{structsubs}

Let $S$ be a subsemigroup of $\AB{n}$ containing 0. Then there is a unique inverse subsemigroup $I \subseteq \AB{n}$ and a unique nilpotent subsemigroup $N \subset \AB{n}$ such that $I$ is a subsemigroup of $S$, $N$ is an ideal of $S$ and $S = I \cup N$.

\et

\proof Define a relation $\equiv_{S}$ on $\{i|(i,i) \in S\}$  by $i\equiv_{S}j$ if and only if both $(i,j),(j,i) \in S$. Since $S$ is a transitive relation by Theorem \ref{subs}, it follows that $\equiv_{S}$ is a partial equivalence relation on $\{1,\ldots,n\}$. By Theorem \ref{subs}, the partition corresponding to $\equiv_{S}$ defines a unique inverse subsemigroup $I$ of $S$ containing $0$.

We claim that $N = (S\backslash I)  \cup {0}$ is a nilpotent ideal of $S$. As a relation, $N$ is irreflexive, since all the non-zero idempotents of $S$ belong to $I$. Let $(i,j) \in N, (k,i) \in S$. If $(k,j) = (k,i)(i,j) \notin N$, then $(k,j) \in I$. Therefore, $(j,k) \in I$ since $I$ is an inverse subsemigroup of $\AB{n}$. Thus, $(j,i) = (j,k)(k,i) \in S$ and $i\equiv_{S}j$. It follows that $(i,j) \in I$ a contradiction. Therefore $N$ is a left ideal of $S$. Dually, it is a right ideal.

Uniqueness follows easily from the fact that all the idempotents of $S$ must belong to $I$.
\qed

Equivalently, Theorem \ref{structsubs} can be given in the language of relations as a description of all transitive relations on an $n$ element set. In this form, it is related to Lemma 1 of \cite{Klaska}. See also \cite{Pfeiffer}. We record this as a corollary.

\bc

Let $R$ be a transitive relation on an $n$ element set. Then there is a unique partial equivalence relation $\Pi$ and a unique  strict partial order $P$ such that $R$ is the disjoint union of $\Pi$ and $P$ and such that $RP \cup PR \subseteq P$.

\ec

\subsection{The BRSC of subsemigroup lattices associated to $\AB{n}$}

\subsubsection{The inverse subsemigroup and nilpotent subsemigroup case}

Theorem \ref{structsubs} motivates looking at the BRSC associated to the subsemigroup lattice of a nilpotent semigroup and the inverse subsemigroup lattice of $\AB{n}$.

We begin with nilpotent semigroups.

\bt \label{nilpbrsc}

Let $N$ be a nilpotent semigroup. Then the subsemigroup complex $\H(N)$ is the uniform matroid, $U_{n,n}$, where $n =|N|$. That is, every subset of $N$ is independent.

\et

\proof Clearly if $S$ is a semigroup of cardinality $n$, then $\H(S)$ is the uniform matroid, $U_{n,n}$, if and only if there exists an ordering of the elements of $S$, $s_{1}, \ldots, s_{n}$ such that the subsemigroup generated by $s_{1}, \ldots, s_{i}$ is just the set $\{s_{1}, \ldots, s_{i}\}$ for $i = 1, \ldots, n$.

We prove by induction on $|N| = n$ that if $N$ is  a nilpotent semigroup then there exists such an ordering of the elements of $N$. If $n \leq 1$, then the assertion is clear. Let $n > 1$. Let $s \in N$ be an element such that $\{s\}$ is a 0-minimal $\mathcal{J}$-class. Then $Ns = sN$ = {0} and it follows that for any nonempty subset, $X \subseteq N$, $\langle X,s\rangle = \langle X\rangle \cup \{s\}$ and that the ideal generated by $\{s\}$ is equal to $\{s,0\}$.

Let $N'$ be the Rees quotient of N by the ideal $\{s,0\}$. By induction, there is an ordering $\{0,s_{1}, \ldots, s_{n-1}\}$ of the elements of $N'$ that proves that $N'$ is an independent set. It follows from the above that $\{0,s,s_{1}, \ldots, s_{n-1}\}$ is an ordering of $N$ that proves that every subset is independent. \qed

We now turn to inverse subsemigroups of $\AB{n}$. More generally, let $S$ be an inverse semigroup. We have defined $ISub(S)$ $(FISub(S))$ to be the lattice of all inverse (full inverse) subsemigroups of $S$. Notice that $FI(S)$ is the interval $[E(S), S]$ of $I(S)$.

These lead us to define two simplicial complexes, $IH(S)$ and $FIH(S)$. Formally, we let $\I_{S}: S \rightarrow ISub(S)$ by letting
$\I_{S}(s)$ be the inverse subsemigroup generated by $s$. This is the same as the subsemigroup of $S$ generated by $\{s,s^{-1}\}$. We have a similar definition in the case of the full inverse subsemigroup lattice by letting $\Phi_{S}: S \rightarrow FISub(S)$ be defined by letting $\Phi_{S}(s)$ be the full inverse subsemigroup generated by $s$, which is the subsemigroup generated by $\{s,s^{-1}\} \cup E(S)$.

We thus get two boolean representable simplicial complexes $IH(S)$ and $FIH(S)$ corresponding to these lattices. As in Section \ref{simply} we can use simplification to identify the atoms of $IH(S)$ $(FIH(S))$ with the set of monogenic (full monogenic) inverse subsemigroups of $S$ and we do so.

In the case of $S = \AB{n}$, $\I(s) =\{s\}$ if $s$ is an idempotent and for $i\neq j$ we get
$$\I((i,j))=\{(i,j),(i,i),(j,i),(j,j),0 \}\thickapprox \AB{2}.$$
We prefer to work with the sublattices $I\sub_{0}(\AB{n})$ and $FI\sub_{0}(\AB{n})$ (which are the intervals from $\{ 0\}$ to the top of the lattices $ISub(S)$ and $FISub(S)$) and their corresponding boolean representable simplicial complexes. By Theorem \ref{subs} and the remark afterward, these are isomorphic to the partition lattices on $n+1$ and $n$ elements respectively. We have the corresponding simplicial complexes $IH_{0}(\AB{n})$ and $FIH_{0}(\AB{n})$. It is well known that the unique simple matroid corresponding to the partial partition lattice on $n$ elements is the graphic matroid on the complete graph $K_n$ (see \cite[Section I.7]{Oxl}). This is the matroid with vertices the edges of $K_n$ and independent sets the forests of $K_n$.

By a minimal generating set $X$ of a semigroup $S$, we mean a generating set such that no proper subset of $X$ generates $S$. When we speak of minimal generating sets of an inverse semigroup, we mean in the variety of inverse semigroups. An oriented spanning tree of a connected graph is a spanning tree along with an orientation on its edges.

\bt \label{mininvgen}

$X$ is a minimal generating set of $\AB{n}$ as an inverse semigroup if and only if the graph $\Gamma(X)$ is an oriented spanning tree of the complete graph $K_n$.

\et

\proof Let $T$ be an oriented spanning tree of $K_n$. If $e$ is an oriented edge of $T$, then we write $e^{-1}$ for the opposite edge. This agrees with the inverse of $e$ thought of as an element of $\AB{n}$. For all $i,j \in \{1, \ldots , n\}$ there is a path from $i$ to $j$. Multiplying the elements on this path (respecting the inverse notation just introduced) in $\AB{n}$ shows that the oriented edges of $T$ generate $\AB{n}$ as an inverse semigroup. If we remove an edge $\{i,j\}$ from $T$, then the resulting graph is not connected and $(i,j)$ is not in the inverse semigroup generated by the remaining edges. Therefore $T$ defines a minimal generating set of $\AB{n}$ as an inverse semigroup.

Conversely, if $X$ is a minimal generating set of $\AB{n}$, then the graph $\Gamma(X)$  corresponding to $X$ must be connected. Furthermore it must be a spanning graph, for if $i$ is not a vertex of $\Gamma(X)$ then no element of the inverse semigroup generated by the elements of $X$ can begin or end in $i$. If $\Gamma(X)$ is not a tree, then there is a cycle from some $i$ to itself in $\Gamma(X)$. Removing an edge from this cycle leaves a connected graph and the remaining edges still generate $\AB{n}$. This contradicts the minimality of $X$. \qed

The arguments in the above proof yield also the following characterization.

\bc

$X$ is a minimal generating set for a (full) inverse subsemigroup of $\AB{n}$ if and only if $\Gamma(X)$ is an oriented (spanning) forest of $K_n$.

\ec

We now prove that the BRSC $IH_{0}(\AB{n})$ is the graphical matroid of $K_n$.

\bt

Let $X$ be a subset of $\AB{n}$. Then $X$ is a face of $IH_{0}(\AB{n})$ if and only if the underlying unoriented graph $\Gamma(X)$ is a forest.

\et

\proof Assume that $\Gamma(X)$ is a forest. Then for any subset $Y$ of $X$, $\Gamma(Y)$ is a forest. Therefore, from Theorem \ref{mininvgen}, any orientation of $\Gamma(Y)$ is a minimal generating set for the inverse semigroup generated by $Y$. It follows that any ordering of the edges of $Y$ gives an ascending chain in the lattice $I\sub_{0}(\AB{n})$ and thus $X$ is a face of $IH_{0}(\AB{n})$.

Assume now that $\Gamma(X)$ is not a forest. Let $e_{1},e_{2}, \ldots , e_{k}$ be any ordering of the edges of $X$. Let $i$ be the maximal index such that the edges
$e_{1}, \ldots , e_{i}$ is a forest, so that $i < k$. By Theorem \ref{mininvgen}, $\{e_{1}, \ldots , e_{i}\}$ is a minimal generating set for the inverse semigroup that it generates. Adding edge $e_{i+1}$ creates a cycle by the definition of $i$. It follows that the inverse subsemigroup generated by $\{e_{1}, \ldots , e_{i}, e_{i+1}\}$ is equal to that of $\{e_{1}, \ldots , e_{i}\}$. Therefore, $X$ is not an independent subset of $IH_{0}(\AB{n})$. \qed

\bc \label{invbrsc}

The BRSC corresponding to $IH_{0}(\AB{n})$ is the graphical matroid of the complete graph $K_n$.

\ec

\subsubsection{The General Case}

We return to the general case of the subsemigroup complex $Sub(\AB{n})$. We first note that the dimension problem was solved in \cite[Theorem 7.1]{CGMP}. In fact, the authors computed the dimension of the subsemigroup lattice of many semigroups including all finite inverse semigroups. We just recall the case of importance for this paper:

\beq
\label{cameron}
\dim \H(\AB{n}) =
\binom{n}{2}+2n-1.
\eeq

We will now characterize the faces of $\H(\AB{n})$ in graph theoretic terms. The following lemma allows us to
omit idempotents from many arguments regarding $\AB{n}$.

\bl
\label{omitid}
Let $X \subseteq \AB{n}$ and let $E = E(\AB{n})$. Then $X \in H(\AB{n})$ if
and only if $X\setminus E \in H(\AB{n})$.
\el

\proof
Since subsets of faces are faces, the direct implication holds
trivially. Assume now that $X\setminus E \in H(\AB{n})$. Then there
exists an enumeration $x_1,\ldots,x_n$ of $X\setminus E$ such that
$$x_1^+ \subset \{ x_1,x_2\}^+ \subset \ldots \subset \{ x_1,\ldots, x_n\}^+$$
holds. Let $e_1,\ldots,e_m$ be an enumeration of $X \cap E$ such that
$e_1 = 0$ if $0 \in X$. It is immediate that
$$e_1^+ \subset \{ e_1,e_2\}^+ \subset \ldots \subset \{ e_1,\ldots,
e_m\}^+.$$
Write $S_i = \{ e_1,\ldots,e_m,x_1, \ldots, x_i\}^+$ for $i =
0,\ldots,n$. Since $se,es \in \{ s,0\}$ for all $s \in \AB{n}$ and $e
\in E$, it follows that $S_i = \{ e_1,\ldots,e_m\}^+ \cup \{ x_1,
\ldots, x_i\}^+$ for each $i$ and so
$$\{ e_1,\ldots,e_m\}^+ = S_0 \subset S_1 \subset \ldots \subset S_n.$$
Hence $X \in H(\AB{n})$ as required.
\qed

Given a directed graph, an edge $a \mapright{} b$ connecting two
distinct vertices is known as
a {\em chord} if there exists some path from $a$ to $b$ avoiding the
edge. Otherwise, it is known as a {\em separating edge} (or
basic edge).

We denote by $K_n$ (respectively $\vec{K}_n$) the complete undirected graph
(respectively the complete directed graph) with vertex set $V_n = \{
1,\ldots,n\}$. Note that we are excluding loops and multiple edges.

Given $X \subseteq \AB{n} \setminus E(\AB{n})$, we define a subgraph
$\Gamma(X)$ of $\vec{K}_n$ with vertex set $V_n$ and edges
$i \mapright{} j$ whenever $(i,j) \in X$.

It is clear that an edge $e$ in a directed graph $G=(V,E)$ is a separating edge if and only if
$e$ is not in the transitive closure of the graph $G\setminus e =(V,E\setminus \{e\})$ considered as
a binary relation. By Theorem \ref{subs}, this is equivalent to saying that $e$ is not
in the subsemigroup generated by $E \setminus \{e\}$, where we consider the edges to be elements of $\AB{n}$.
This suggests that there is a close connection between separating edges and faces in $\H(\AB{n})$. We make this
precise in the next lemma.

\bl
\label{facesh}
The following conditions are equivalent for a given $X \subseteq
\AB{n} \setminus E(\AB{n})$:
\bi
\item[(i)] $X \in H(\AB{n})$;
\item[(ii)] there exists an enumeration $e_1,\ldots, e_m$ of the edges
  of $\Gamma(X)$ such that, for $i = 1,\ldots,m$, $e_i$ is a
  separating edge of the subgraph of $\Gamma(X)$ obtained by removing
  $e_{i+1},\ldots,e_m$;
\item[(iii)] every non edgeless subgraph of $\Gamma(X)$ has a
  separating edge.
\ei
\el

\proof
(i) $\Rw$ (iii). Let $x_1,\ldots,x_n$ be an enumeration of $X$ such that
\beq
\label{facesh1}
x_1^+ \subset \{ x_1,x_2\}^+ \subset \ldots \subset \{ x_1,\ldots,
x_n\}^+.
\eeq

Let $\Gamma' = (V',E')$ be a subgraph of $\Gamma(X)$ with $E' \neq
\emptyset$. Let $x_{i_1},\ldots,x_{i_m}$ be the subsequence of
$x_1,\ldots,x_n$ corresponding to the edges of $X'$, and write
$x_{i_m} = (a,b)$. We claim that $a \mapright{} b$ is a separating
edge of $\Gamma'$.

Indeed, suppose $a \mapright{} b$ is a chord of $\Gamma'$. Then there
exists an alternative path
$$a = c_0 \mapright{} c_1 \mapright{} \ldots \mapright{} c_k = b$$
in $\Gamma'$ avoiding the edge $a \mapright{} b$. It follows that
$$(c_0,c_1), \ldots, (c_{k-1},c_k) \in \{
x_{i_1},\ldots,x_{i_{m-1}} \}^+$$
and so
$$x_{i_m} = (a,b) = (c_0,c_1) \ldots (c_{k-1},c_k) \in \{
x_{i_1},\ldots,x_{i_{m-1}} \}^+,$$
contradicting (\ref{facesh1}). Therefore $a \mapright{} b$ is a separating
edge of $\Gamma'$.

(iii) $\Rw$ (ii). By successive application of condition (iii).

(ii) $\Rw$ (i). Let $e_1,\ldots, e_m$ be an enumeration of the edges
  of $\Gamma(X)$ such that, for $i = 1,\ldots,m$, $e_i$ is a
  separating edge of the subgraph of $\Gamma(X)$ obtained by removing
  $e_{i+1},\ldots,e_m$. Let $e_i$ be the edge $a_i \mapright{} b_i$
  and write $x_i = (a_i,b_i)$. Then $x_1,\ldots, x_m$ is an
  enumeration of the elements of $X$. Clearly, $x_i \in \{
  x_1,\ldots,x_{i-1}\}^+$ implies that there exists a path from $a_i$
  to $b_i$ using only the edges $e_1,\ldots,e_{i-1}$, a
  contradiction. Therefore (\ref{facesh1}) holds and so $X \in
  H(\AB{n})$.
\qed

The concept of separating edge allows us to develop a recursive
procedure to construct the faces of $\H(\AB{n})$. We first give a graph theoretic formulation and then reformulate it in terms of the connections to subsemigroups of $\AB{n}$ from Theorem \ref{subs}.

\bl
\label{recur}
Let $V_n = W_1 \cup W_2$ be a nontrivial partition. Let $X_i
\subseteq \AB{n} \setminus E(\AB{n})$ be a face of $\H(\AB{n})|_{W_i}$ and
let $a_i \in W_i$ for $i = 1,2$. Let
\beq
\label{recur1}
X = X_1 \cup X_2 \cup \{ (a_1,a_2) \} \cup Y
\eeq
with $Y \subseteq W_2 \times W_1$. Then:
\bi
\item[(i)] $X$ is a face of $H(\AB{n})$;
\item[(ii)] every nonempty face of $H(\AB{n})$ containing no idempotents
  can be obtained this way.
\ei
\el

\proof
(i) Note that $X \subseteq \AB{n} \setminus E(\AB{n})$ and $a_1
\mapright{} a_2$ is a separating edge of $\Gamma(X)$. Let $X' = X
\setminus \{ (a_1,a_2)\}$. By Lemma
\ref{facesh}, it suffices to show that every non edgeless subgraph of
$\Gamma(X')$ has a separating edge. Let $\Gamma'$ be such a subgraph.

Assume first that $\Gamma'$ has edges with both endpoints in
$W_1$. Let $\Gamma'_1$ be the subgraph of $\Gamma'$ induced by all
such edges. Since $\Gamma'_1$ is also a non edgeless subgraph of
$\Gamma(X_1)$ and $X_1$ is a face, it follows from Lemma \ref{facesh}
that $\Gamma'_1$ has
a separating edge, which must also be a separating edge of $\Gamma'$.

Similarly, we deal with the case where $\Gamma'$ has edges with both
endpoints in $W_2$. Thus we may assume that all edges of $\Gamma'$ are
contained in $W_2 \times W_1$, but then every edge is trivially a
separating edge.

(ii) Let $X \subseteq \AB{n} \setminus E(\AB{n})$ be a nonempty face. By
Lemma \ref{facesh}, $\Gamma(X)$ has a separating edge $a_1 \arr
a_2$. Let $X' = X \setminus \{ (a_1,a_2) \}$ and let $W_1$ denote the
set of all vertices $b \in V_n$ such that there exists a path $a_1 \arr b$ in $\Gamma(X')$
(including the trivial path so $a_{1} \in W_1$). Finally, let $W_2 = V_n \setminus
W_1$. Since $a_i \in W_i$ for $i = 1,2$, $W_1 \cup W_2$ is a
nontrivial partition of $V_n$. For $i = 1,2$, let
$$X_i = \{ (p,q) \in X \mid p,q \in W_i\}.$$
Then $X_i \subseteq \AB{n} \setminus E(\AB{n})$ is a face of
$\H(\AB{n})|_{W_i}$. We show that (\ref{recur1}) holds for
$Y = X \setminus \{ X_1 \cup X_2 \cup \{ (a_1,a_2) \} \}$.

Indeed, suppose that $(b_1,b_2) \in X$ and $b_i \in W_i$. Then there
exists a path $a_1 \arr b_1$ and therefore a path $a_1 \arr
b_2$ in $\Gamma(X)$. Since $b_2 \notin W_1$, then there is no path $a_1 \arr
b_2$ in $\Gamma(X')$ and so
we must have $(b_1,b_2) = (a_1,a_2)$.
This completes the proof of the claim and therefore of
the lemma.
\qed

We can now settle problem (P5). Simultaneously, we provide an
alternative proof for the aperiodic case in (\ref{cameron}).

\bp
\label{pure}
The complex $\H(\AB{n})$ is pure of dimension $\binom{n}{2}+2n-1$
for every $n \geq 1$.
\ep

\proof
We proceed by induction on $n$. The case $n = 1$ being trivial, assume
that $n > 1$ and the claim holds for smaller values. Let $F$ be a
facet of $\H(\AB{n})$. By Lemma \ref{omitid}, we must have $E(\AB{n})
\subseteq F$. Let $X = F \setminus E(\AB{n})$.
By Lemma \ref{recur}, there exists a nontrivial partition $V_n = W_1
\cup W_2$ and a (disjoint) decomposition (\ref{recur1}) such that $X_i \subseteq
W_i$ for $i = 1,2$ and $Y \subseteq
W_2 \times  W_1$.

Now
\bi
\item
$X_i$ must be maximal among the faces contained in $(W_i \times
  W_i) \setminus E(\AB{n})$,
\item
$Y = W_2 \times  W_1$,
\ei
otherwise Lemma \ref{recur} would grant us a face of $\AB{n}$ strictly
containing $F$.

Write $m_i = |W_i|$. By Lemma \ref{omitid} and the induction hypothesis,
$X_i \cup \{ (a,a) \mid a \in W_i\} \cup \{ 0 \}$
is a face of dimension $\binom{m_i}{2} + 2m_i-1$, and so
$|X_i| = \binom{m_i}{2} + m_i-1$. Now (\ref{recur1})
yields
$$|F| = |X| + n+1 = \binom{m_1}{2} + m_1-1 + \binom{m_2}{2} + m_2-1 +
1 + m_2m_1 +n+1.$$
Since $m_2 = n-m_1$, a straightforward computation yields $|F| =
\binom{n}{2} + 2n$ and we are done.
\qed

We can now derive some information on the facets of $\AB{n}$. The next lemma is the graph theoretic interpretation of the construction of maximal subsemigroups for arbitrary finite 0-simple semigroups \cite{Graham, GGR} in the special case of aperiodic Brandt semigroups. We
introduce some more graph-theoretical tools.

Let $\Gamma = (V,E)$ be an undirected graph. An {\em orientation} of
$\Gamma$ is a binary relation $O \subseteq V \times V$ such that
$$\begin{array}{rcl}
O&\to&E\\
(p,q)&\mapsto&\{p,q\}
\end{array}$$
is a bijection. Intuitively, we are choosing, for each undirected edge
$p \edge q$, one of the directed edges $p \mapright{} q$, $q
\mapright{} p$. An orientation is {\em acyclic} if it contains no
directed cycle. An acyclic orientation determines a partial order $\leq_O$ on $V$
by declaring $v \leq_{O} w$ if there is a path (including the empty path) from $v$ to $w$ in the directed graph $O$.
Conversely, the directed graph of a partial order $\leq$ is an acyclic orientation, which is just the transitive closure
of the Hasse Diagram of $\leq$.

Consider the complete undirected graph $K_n$. A {\em spanning tree} of
$K_n$ is a subtree $T$ containing all $n$ vertices (and having therefore
$n-1$ edges). We denote by $K_n\setminus T$ the graph obtained by
removing from $K_n$ all the edges of $T$.

\bp
\label{facets}
Every facet of $\H(\AB{n})$ is of the form
\beq
\label{facets2}
E(\AB{n}) \cup \{ (p,q) \mid \{p,q\} \in T\} \cup
\{ (p,q) \mid (p,q) \in O \},
\eeq
where $T$ is a spanning tree of $K_n$ and $O$ is an acyclic orientation
of $K_n\setminus T$.
\ep

\proof
Given $p,q \in V_n$ distinct and $X \subseteq \AB{n} \setminus E(\AB{n})$,
we say that $p \edge q$ is a {\em link} of $\Gamma(X)$ if both $p
\mapright{} q$ and $q \mapright{} p$ are edges of $\Gamma(X)$.

We use induction on $n$. The case $n = 1$ being trivial, assume that $n > 1$, $F$ is a facet of $\H(\AB{n})$ and the claim holds for $m < n$. By the proof of Proposition \ref{pure}, there exists a nontrivial partition
$V_n = W_1 \cup W_2$, $a_i \in W_i$ and a facet $F_i$ of $\H(\AB{n})|_{W_i}$ for $i = 1,2$ such that
$$F = E(\AB{n}) \cup F_1 \cup F_2 \cup \{ (a_1,a_2) \} \cup (W_2 \times W_1).$$
By the induction hypothesis, we may write
$$F_i = E(\AB{W_i}) \cup \{ (p,q) \mid \{p,q\} \in T_i\} \cup
\{ (p,q) \mid (p,q) \in O_i \},$$
where $T_i$ is a spanning tree of $K_{W_i}$ and $O_i$ is an acyclic orientation
of $K_{W_i}\setminus T_i$ for $i = 1,2$.

Since $(a_2,a_1) \in F$ as well, we define
$$T = T_1 \cup T_2 \cup \{ \{ a_1,a_2\} \},$$
which is obviously a spanning tree for $K_n$. It is easy to check that
$$O = O_1 \cup O_2 \cup (W_2 \times W_1)$$
is an acyclic orientation of $K_n \setminus T$. Therefore $F$ equals (\ref{facets2}) and we are done.
\qed

The converse of Proposition \ref{facets} does not
hold. The smallest counterexample is given by the following graph on $\{1,2,3,4\}$. One easily checks that this graph does not represent a face and thus not a facet  of $\mathcal{H}\AB{4})$. To make pictures easier to interpret, we shall represent an oriented
edge $p \arr q$ in black, and a link $p \edge q$ in red.

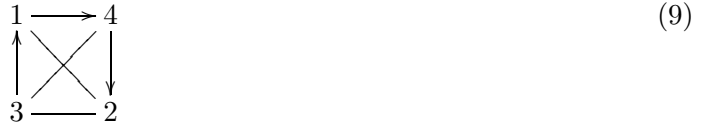
\beq
\label{square}
\xymatrix{
1 \ar[r] & 4 \ar[d] \\
3 \ar[u] \ar@{-}@[red][ur] & 2 \ar@{-}@[red][l] \ar@{-}@[red][ul]
}
\eeq

We present next a straightforward construction of a facet: we take
links
$$1 \edge 2 \edge \ldots \edge n$$
plus all edges of the form $i \arr j$ with $i > j+1$. The conditions
of Lemma \ref{recur} are satisfied by the partition $V_n = V_{n-1}
\cup \{ n \}$, so it follows easily by induction that this graph
defines a face $X_n$ of $\H(\AB{n})$. Hence $X_n \cup E(\AB{n})$ is a
face. Since $|X_n \cup E(\AB{n})| = \binom{n}{2} + 2n$, it follows from
(\ref{cameron}) that $X_n \cup E(\AB{n})$ is indeed a facet of $\H(\AB{n})$.

We settle next problem (P6).

\bp
\label{matro}
The complex $\H(\AB{n})$ is a matroid if and only if $n \leq 3$.
\ep

\proof
The case $n \leq 2$ is trivial since by (\ref{cameron}) we have
$\dim \H(\AB{n}) = \binom{n}{2} +2n-1 = n^2 = |\AB{n}|-1$ and so every
subset of $\AB{n}$ must be a face.

Assume now that $n = 3$. Then the facets of $\AB{3}$ correspond to
graphs of the form
$$
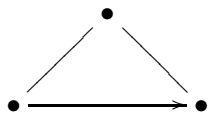
\xymatrix{
& \bullet \ar@{-}@[red][dl] \ar@{-}@[red][dr] & \\
\bullet \ar[rr] && \bullet
}$$
Let
$$Z = \{ (i,j) \mid i,j \in V_3,\; i \neq j\}.$$
It follows from Proposition \ref{pure} that
\beq
\label{matro1}
H(\AB{3}) = \{ X \subseteq \AB{3} \mid Z \not\subseteq X \}.
\eeq

Assume that $I,J \in H(\AB{3})$ are such that $|I| = |J|+1$. We may assume that $J \cap E(\AB{n}) \supseteq I \cap E(\AB{n})$, otherwise we may add some idempotent in $I \setminus J$ to $J$ and still get a face.

If $|J \cap Z| =
5$, then $J \cap E(\AB{n}) \supseteq I \cap E(\AB{n})$ implies $Z \subseteq I$, contradicting (\ref{matro1}). Hence
$|J \cap
Z| < 5$ and so $J \cup \{ p \} \in H(\AB{3})$
for every $p \in I \setminus J$ in
view of (\ref{matro1}). Therefore $\H(\AB{3})$ is a matroid.

Assume now that $n \geq 4$. Let $I,J \subseteq \AB{n} \setminus E(\AB{n})$
be defined by the graphs
$$
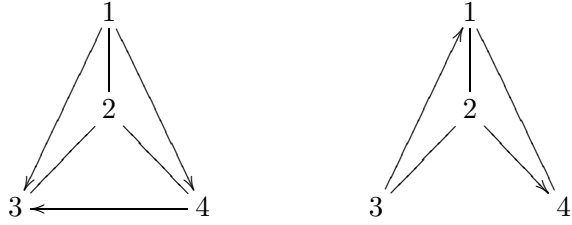
\xymatrix{
&1 \ar[ddl] \ar[ddr] \ar@{-}@[red][d] &
&&
&1 \ar@{-}@[red][ddr] \ar@{-}@[red][d] &\\
&2 \ar@{-}@[red][dl] \ar@{-}@[red][dr] &
&&
&2 \ar@{-}@[red][dl] \ar[dr] &\\
3 && 4 \ar[ll]
&&
3 \ar[uur] && 4
}$$
respectively.

Considering the partitions $V_n = \{ 3 \} \cup (V_n \setminus \{ 3\})$
and $V_n = \{ 4 \} \cup (V_n \setminus \{ 4\})$, respectively, we
deduce from Lemma \ref{recur} that $I,J \in H(\AB{n})$. Clearly, $|I| =
|J|+1$. Now $|I \setminus J| = 3$ and $J \cup \{ i \}$, for $i \in I
\setminus J$, produces the three graphs
$$
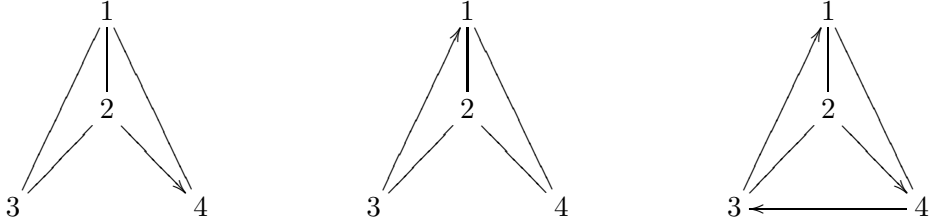
\xymatrix{
&1 \ar@{-}@[red][ddr] \ar@{-}@[red][d] &
&&
&1 \ar@{-}@[red][ddr] \ar@{-}@[red][d] &
&&
&1 \ar@{-}@[red][ddr] \ar@{-}@[red][d] &\\
&2 \ar@{-}@[red][dl] \ar[dr] &
&&
&2 \ar@{-}@[red][dl] \ar@{-}@[red][dr] &
&&
&2 \ar@{-}@[red][dl] \ar[dr] &\\
3 \ar@{-}@[red][uur] && 4
&&
3 \ar[uur] && 4
&&
3 \ar[uur] && 4 \ar[ll]
}$$
In view of Proposition \ref{facets}, the links of a face must constitute a forest. Thus $J \cup \{ i \}$ is not a face in the two first
cases. On the other hand, the third graph is (\ref{square}), and we
have already established that it does not correspond to a
face. Therefore $\H(\AB{n})$ is not a matroid if $n \geq 4$.
\qed

We now present an alternative approach for constructing faces and facets in $\mathcal{H}(\AB{n})$. This is based on Theorem \ref{subs} and based on a global understanding of the inductive approach inspired by \cite{GGR} and that appeared in the graph theoretic approach above in Theorem \ref{recur}. This approach gives a proof that a set of elements in $\AB{n}$ is a face, by listing an order in which the elements of the set are a transversal of successive differences. At the end of this discussion, we will give the connection between the algebraic and the graph theoretic approaches.

Let $L$ be a strict linear order of the set $\{1, \ldots, n\}$ and let $T$ be a spanning tree of $K_n$. Then there is a unique orientation $\mathcal{O}$ of the edges of $T$ that is opposite to the order given by $L$. We define a facet $\Phi(L,T,\mathcal{O})$ as follows. By the remarks after Theorem \ref{nilpsub}, $L$ corresponds to a unique maximal nilpotent subsemigroup $N(L)$ of $\AB{n}$. Clearly, $|N(L)| = {n\choose 2}+1$. By Theorem \ref{nilpbrsc} $N(L)$ is a face of $\mathcal{H}(\AB{n})$. We can now add the $n$ non-zero idempotents $(i,i), i=1, \ldots, n$ one at a time to get a face of size ${n\choose 2}+n+1$, that defines a chain in $\text{Sub}(\AB{n})$ whose top is the subsemigroup $N(L) \cup E(\AB{n})$. Finally, the oriented spanning tree $(T,\mathcal{O})$ defines a unique minimal generating set $X(T,\mathcal{O})$ of $\AB{n}$ as an inverse semigroup by
Theorem \ref{mininvgen} and $X(T,\mathcal{O})$ has size $n-1$. We define the set $\Phi(L,T,\mathcal{O})=N(L) \cup E(\AB{n}) \cup X(T,\mathcal{O})$ of size ${n\choose 2}+2n$. It follows from \cite{CGMP}, or from Proposition \ref{pure} that this is the size of a maximal possible face.

Note that $\mathcal{O}$ is fully determined by $L$ and $T$, so $\Phi(L,T,\mathcal{O})$  ia actually a function of $L$ and $T$, but we opt for our notation for emphasis.

We give an example of this construction.

\be

Let $n=4$. Let $L$ be the linear order $ 2 < 4 < 3 < 1$ and let $T$ be the spanning tree of $K_4$ with edges $\{\{3,2\},\{2,1\},\{1,4\}\}$. Then
$\mathcal{O}=\{(3,2),(1,2),(1,4)\}$. The set $\Phi(L,T,\mathcal{O})$ is equal to $\{0,(2,1), (2,3), (4,1), (2,4), (4,3), (3,1), (1,1), (2,2), (3,3), (4,4), (3,2), (1,4), (1,2)\}$. The reader can verify that adding elements one at a time from left to right gives a chain of subsemigroups that proves that this set is a facet.

\ee

However, not every set formed this way is a face as the following example shows.

\be

Let $n=4$. Let $T$ be the spanning tree of $K_4$ with edges $\{(1,2),(2,3),(3,4)\}$ and $L$ the linear order $3 < 1 < 4 < 2$. Then $\mathcal{O} =\{(2,1),(2,3),(4,3)\}$. One computes that $\Phi(L,T,\mathcal{O}) \setminus E(\AB{n})=\{(3,2),(3,4),(1,2),(3,1),(1,4),(4,2),(2,1),(2,3),(4,3)\}$ and these are precisely the edges of the graph (\ref{square}). Thus $\Phi(L,T,\mathcal{O})$ is not a face of $\mathcal{H}(\AB{n})$.

\ee

Despite this we now prove that every facet of the subsemigroup complex of $\AB{n}$ is of the form $\Phi(L,T,\mathcal{O})$. Then we explore the connection between this construction and the graph theoretic construction described previously. The proof below is part of Proposition 3.1 of \cite{CGMP}.

\bl \label{butter}

Let $U,V$ be subsemigroups of a finite semigroup $S$ such that $U$ is covered by $V$ in the subsemigroup lattice of $S$. Let $I$ be an ideal of $S$. Then either
$U \cap I = V \cap I$ and $U\backslash I$ is covered by $V\backslash I$ in the subsemigroup lattice of the Rees quotient $S\backslash I$ or $U\backslash I = V\backslash I$ and $U \cap I$ is covered by  $V \cap I$ in the subsemigroup lattice of $I$.

\el

\proof In any semigroup, the union of a subsemigroup and an ideal is a subsemigroup. Thus, if $U \leq V$ in the semigroup lattice of $S$ and $I$ is an ideal of $S$, then $(V \cap I)$ is an ideal of $V$ and we
have the inequalities, $ U \leq U \cup (V \cap I) \leq V$ in the subsemigroup lattice of $S$. Since $V$ covers $U$, we have either that $U = U \cup (V \cap I)$ or that $U \cup (V \cap I) = V$. In the first case we have $U \cap I = V \cap I$ and $U\backslash I$ is covered by $V\backslash I$ in the subsemigroup lattice of the Rees quotient $S\backslash I$. In the second case we have that $U\backslash I = V\backslash I$ and $U \cap I$ is covered by  $V \cap I$ in the subsemigroup lattice of $I$. \qed

This allows us to build facets of the subsemigroup complex of a semigroup $S$ from facets of an ideal $I$ and those of the Rees quotient $S \backslash I$.

\bt \label{ideal}

Let $I$ be an ideal of a finite semigroup $S$. Let $F_1$ be a facet of the subsemigroup complex of $I$ and $F_2$ be a facet of the subsemigroup complex of $S \backslash I$. Then $F_1 \cup F_2$ is a facet of the subsemigroup complex of $S$. Conversely, if $F$ is a facet of the subsemigroup complex of $S$, then $F \cap I$ is a facet of the subsemigroup complex of $I$ and $F \backslash I$ is a facet of the subsemigroup complex of $S \backslash I$.

\et

\proof Clearly if $F_1$ is a facet of the subsemigroup complex of $I$ and $F_2$ is a facet of the subsemigroup complex of $S \backslash I$, then $F_1 \cup F_2$ is a facet of the subsemigroup complex of $S$.

Conversely, let $F$ be a facet of the subsemigroup complex of $S$. Then there is an ordering $s_{1}, \ldots , s_{k}$ of the elements of $F$ such that
$S_{0} = \emptyset < S_{1} < \ldots < S_{i} < S_{i+1} < \ldots < S_{k} = S$ is a maximal chain in the subsemigroup lattice of $S$, where $S_{i}$ is the subsemigroup generated by $\{s_{1}, \ldots s_{i}\}$. Therefore, $S_{i+1}$ covers $S_{i}$ for $i=0, \ldots , k-1$. By Theorem \ref{butter}, for each $0 \leq i \leq k-1$, either $S_{i+1} \backslash S_{i} \subseteq I$ or $S_{i+1} \backslash S_{i} \subseteq S \backslash I$. In the first case, $s_{i+1} \in I$ and $S_{i+1} \cap I$ covers $S_{i} \cap I$ in the semigroup lattice of $I$. In the second case, $s_{i+1} \in S \backslash I$ and $S_{i+1} \backslash I$ covers $S_{i} \backslash I$ in the semigroup lattice of $S \backslash I$. It follows that $F \cap I$ is a facet in the subsemigroup complex of $I$ and $F \backslash I$ is a facet in the subsemigroup complex of $S \backslash I$. \qed

\bt \label{sgpfacets}

Let $F$ be a subset of $\AB{n}$. If $F$ is a facet of the subsemigroup complex of $\AB{n}$ then there exists a linear order $L$, a spanning tree $T$ of $K_n$ and an orientation $\mathcal{O}$ of $T$ such that $F=\Phi(L,T,\mathcal{O})$

\et

\proof
Let $F$ be a facet of the subsemigroup complex of $\AB{n}$. Then there is an ordering $s_{1}, \ldots , s_{k}$ of the elements of $F$ such that
$S_{0} = \emptyset < S_{1} < \ldots < S_{i} < S_{i+1} < \ldots < S_{k} = \AB{n}$ is a maximal chain in the subsemigroup lattice of $\AB{n}$, where $S_{i}$ is the subsemigroup generated by $\{s_{1}, \ldots s_{i}\}$. Therefore, $S_{i+1}$ covers $S_{i}$ for $i=0, \ldots , k-1$.

It follows that $S_{k-1}$ is a maximal subsemigroup of $\AB{n}$ and that $F' = F \backslash \{s_{k}\}$ is a facet of $S_{k-1}$. By \cite{GGR}, there is a non-trivial partition $X,Y$ of $\{1, \ldots n\}$ such that $S_{k-1} = \AB{X} \cup \AB{Y} \cup (Y \times X)$ and $s_k \in X \times Y$. Clearly, $I = (Y \times X) \cup \{0\}$ is a nilpotent (in fact a 0-semigroup) ideal of $S_{k-1}$. It follows from Theorem \ref{ideal} that $F \cap I$ is a facet of $I$ and that $F \backslash I$ is a facet of $S \backslash I$.

By Theorem \ref{nilpbrsc}, $(F \cap I) = I$. The graph $\Gamma(I)$ corresponding to $I$ is the complete bipartite tournament $T(Y,X)$, that is, the directed graph whose edges are all the edges from every vertex in $Y$ to every vertex in $X$.

The Rees quotient $S \backslash I$ is the 0-disjoint union of $\AB{X}$ and $\AB{Y}$. Since $\AB{Y}$ is an ideal in $S \backslash I$, we can appeal again to Theorem \ref{ideal} to write $F \backslash I = (F_{X} \cup F_{Y})$, where $F_{X}$ ($F_{Y}$) is a facet of $\AB{X}$ ($\AB{Y}$). By induction, $F_{X}= \Phi(L_{X},T_{X},\mathcal{O}_{X})$ ($F_{Y}= \Phi(L_{Y},T_{Y}, \mathcal{O}_{Y})$) for appropriate linear orders, spanning trees and orientations for $X$ and $Y$ appropriately. Then,
$L = L_{X} \cup L_{Y} \cup (Y \times X)$
is a linear order of $\{1, \ldots , n\}$. If $s_k = (a,b)$, then $T=T_{X} \cup T_{Y} \cup \{ \{ a,b\} \}$ is a spanning tree of $\{1, \ldots , n\}$ with orientation opposite to $L$,  $\mathcal{O} = \mathcal{O}_{X} \cup \mathcal{O}_{Y} \cup \{ (a,b) \}$ and $F = \Phi(L,T,\mathcal{O})$.
\qed

We now give the connection between the semigroup construction of facets of $\mathcal{H}(\AB{n})$ and the graph theoretic approach. Suppose that a facet $F$ is defined as in Theorem \ref{sgpfacets} as $F=\Phi(L,T,\mathcal{O})$. For each edge $(i,j) \in \mathcal{O}$, the edge $(j,i)\in L$ by construction. Thus, $L \setminus \{(j,i)|(i,j) \in \mathcal{O}\}$ is a poset on $\{1,\ldots n\}$ and corresponds to a unique acyclic orientation of $K_{n} \setminus T$. Clearly, the edges in the graph whose links are $T$ and this oriented acyclic graph is $F \setminus E(\AB{n})$. Conversely, let a facet $F$ be such that $F\setminus E(\AB{n})$ is given by a spanning tree $T$ of $K_n$ and an oriented acyclic graph with edges $K_{n} \setminus T$. Let $L$ be any extension of the poset corresponding to this orientation of $K_{n} \setminus T$ to a linear order. We define $\mathcal{O}$ to be the orientation of $T$ that is opposite to $L$. Then $F$ is defined by $F=\Phi(L,T,\mathcal{O})$. We record this in the following Theorem.

\bt \label{connection}

Let $F$ be a facet and let $F\setminus E(\AB{n})$ be given by a tree of links $T$ and an acyclic orientation on $K_{n} \setminus T$. Then the number of representations of $F$ in the form $F=\Phi(L,T,\mathcal{O})$ is equal to the number of extensions of the poset on $\{1,\ldots n\}$ given by the acyclic orientation on $K_{n} \setminus T$ to a linear order on $\{1,\ldots n\}$.

\et

The number of linear extensions of a poset on $\{1,\ldots n\}$ is a well studied and difficult to understand combinatorial number. In Subsection \ref{asy} we give an estimate on the number of facets in $\mathcal{H}(\AB{n}$.

\subsection{Counting facets: asymptotics}
\label{asy}

It seems out of reach to count the exact number of facets
of $\H(\AB{n})$ for arbitrary $n$. We recall that the lattice of subsemigroups of $B(n)$ containing 0 is isomorphic to the lattice of transitive relations on an $n$-element set. There is no closed formula for the number of such relations and thus we don't even have a formula for the size of the subsemigroup lattice of $B(n)$.

It is nevertheless possible to get some results of asymptotic type.
Let $f_n$ denote the number of facets of $\H(B(n))$.

\bt
\label{countf}
\bi
\item[(i)]
$f_n < n!n^{n-2}$ for every $n \geq 2$;
\item[(ii)]
$f_n > 2^{n-2}n^{n-2}$ for every $n \geq 4$.
\ei
\et

\proof
(i) By Theorem \ref{sgpfacets}, we have $f_n \leq t_n\ell_n$, where $t_n$ is the number of spanning trees of $K_n$ and $\ell_n$ is the number of linear orderings of $\{ 1,\ldots,n\}$. Now $t_n = n^{n-2}$
by Cayley's Theorem (see \cite[Theorem 5.2.3]{CM}) and $\ell_n = n!$. Hence $t_n \leq n!n^{n-2}$.

Suppose that $n \geq 2$ and $T$ is the (unique) spanning tree of $K_n$ where $n$ has degree $n-1$. Then the position of $n$ in the linear ordering is irrelevant and it is easy to check (out of symmetry) that the facets with spanning tree $T$ correspond precisely to the linear orderings of $\{ 1,\ldots,n-1\}$, hence there are precisely $(n-1)!$ of them. This implies in particular that $f_n < n!n^{n-2}$ for $n \geq 2$.

%
%

(ii) We show that, for a given spanning tree $T$ of $K_n$, there exist at least $2^{n-2}$ facets of $\H(B(n))$ containing the spanning tree $T$ (see Proposition \ref{facets}).

An {\em exfoliation} of $T$ is an enumeration $X = (i_1,\ldots,i_n)$ of the elements of $\{ 1,\ldots,n\}$ such that $i_j$ is a leaf of the subtree of $T$ induced by $\{ i_1,\ldots,i_j\}$ for $j = n,\ldots,2,1$. Note that by erasing a leaf from a tree we always get a tree. Since every tree with more than one vertex has at least two leaves, it follows that $T$ admits at least $2^{n-1}$ exfoliations (exactly $2^{n-1}$ if $T$ is a linear graph).

Given an exfoliation $X = (i_1,\ldots,i_n)$  of $T$, we define
$$F(X) = E(B(n)) \cup \{ (p,q) \mid \{p,q\} \in T\} \cup
\{ (i_j,i_k) \mid  1 \leq j < k \leq n \}.$$
We claim that $F(X)$ is a facet of $\H(B(n))$. Note that it is obvious that the latter set defines an acyclic orientation of $K_n \setminus T$, but we have already remarked that this is not enough to produce a face.

Since $|F(X)| = \binom{n}{2}+2n-1 = \dim\H(B(n))$ by Proposition \ref{pure}, it is enough to show that $F(X) \in H(B(n))$. By Lemma \ref{omitid}, we only need to show that
\beq
\label{countf1}
\{ (p,q) \mid \{p,q\} \in T\} \cup
\{ (i_j,i_k) \mid  1 \leq j < k \leq n \} \in H(B(n)).
\eeq

For $r = 1,\ldots,n$, let $X_r = \{ i_1,\ldots,i_r\}$ and
$$Y_r = \{ (i_j,i_k) \mid  1 \leq j < k \leq r \} \cup \{ (i_k,i_j) \mid  1 \leq j < k \leq r \mbox{ and } \{ i_j,i_k\} \in T \}.$$
We show that $Y_r \in H(B(n))$ by induction on $r$. Since $Y_1 = \emptyset$, the claim holds trivially for $r = 1$.

Assume now that $1 < r \leq n$ and $Y_{r-1} \in H(B(n))$.
Consider the partition $X_r = \{ i_r \} \cup X_{r-1}$. Trivially, $\emptyset$ is a face of $\H(B(n))|_{\{ i_r \}}.$ On the other hand, it follows from the induction hypothesis that
$Y_{r-1}$ is a face of $\H(B(n))|_{X_{r-1}}.$
Since $i_r$ is a leaf of the subtree of $T$ induced by $\{ i_1,\ldots,i_r\}$, there exists a unique $j \in \{ 1,\ldots,r-1\}$ such that $i_j \edge i_r$ is an edge of $T$. Thus we may write
$$Y_r = \emptyset \cup Y_{r-1} \cup \{ (i_r,i_j) \} \cup \{ (i_k,i_r) \mid 1 \leq k < r \},$$
and it follows from Lemma \ref{recur} that $Y_r \in H(B(n))$.

In particular, $F(X) \setminus E(B(n)) = Y_n \in H(B(n))$ and so (\ref{countf1}) holds. Therefore $F(X)$ is a facet of $\H(B(n))$.

Let $X = (i_1,\ldots,i_n)$ and $X' = (i'_1,\ldots,i'_n)$ be exfoliations of $T$ such that $i_s \neq i'_s$ for some $s \geq 3$. Without loss of generality, we may assume that $s$ is maximum. We claim that $F(X) \neq F(X')$.

Indeed, $(i'_s,i_s) \in F(X)$ since $i'_s$ must necessarily occur before $i_s$ in $X$. Suppose that $(i'_s,i_s) \in F(X')$. Since $i_s$ occurs before $i_{s'}$ in $X'$, it follows that $i_s \edge i'_s$ is an edge of $T$. Since $i_s$ and $i'_s$ are supposed to be both leaves of the subtree of $T$ induced by $\{ i_1,\ldots,i_s \} = \{ i'_1,\ldots,i'_s \}$, it follows that $s = 2$, a contradiction. Thus $(i'_s,i_s) \notin F(X')$ and so $F(X) \neq F(X')$.

We had already remarked that $T$ admits at least $2^{n-1}$ exfoliations, and only two exfoliations may share the last $n-2$ components. Thus we have at least $2^{n-2}$ facets in $\H(B(n))$ containing the spanning tree $T$, and so $f_n \geq 2^{n-2}n^{n-2}$ in view of Cayley's Theorem.

Since $n \geq 4$, we have $2^{n-2} < (n-1)!$. Since we showed in the proof of part (i) that there exist $(n-1)!$ facets in $\H(B(n))$ containing a spanning tree with a vertex of degree $n-1$, we get $f_n > 2^{n-2}n^{n-2}$.
\qed

\bc
\label{logfac}
\bi
\item[(i)] $\frac{\log f_n}{n\log n} < 2$ for every $n \geq 2$;
\item[(ii)] $\frac{\log f_n}{n\log n} > 1$ for every $n \geq 8$.
\ei
\ec

\proof
(i) Since $n! < n^n$ for every $n \geq 2$, it follows from Theorem
\ref{countf}(i) that $f_n < n^{2n-2} < n^{2n}$. Hence $\log f_n < \log
(n^{2n}) = 2n\log n$ and so $\frac{\log f_n}{n\log n} < 2$.

(ii) By Theorem \ref{countf}(ii), we have $f_n > 2^{n-2}n^{n-2}$ for $n \geq 4$. Since
$2^6 = 8^2$, we have $2^{n-2} \geq n^2$ for every $n \geq 8$, yelding
$f_n > n^{n}$. Thus $\log f_n > \log
(n^{n}) = n\log n$ and so $\frac{\log f_n}{n\log n} > 1$ for every $n
\geq 8$.
\qed

\section{Suggested future problems}

The work in this paper can of course be carried out for the subsemigroup complex of an arbitrary finite semigroup. Another important project would be a similar analysis of the BRSC corresponding to the congruence lattice of a finite semigroup. The work can be further generalized by a study of the BRSC corresponding to subalgebra and congruence lattices of an arbitrary finite Universal Algebra.

More specifically, it would be very useful to look at the problems studied in this paper for $B(1,n)$ for an arbitrary finite Brandt semigroup $B(G,n)$ over a finite group $G$. In our paper \cite{MRS2} we show that the BRSC defined by the Rhodes lattice, which is the semilattice of aperiodic inverse subsemigroups of $B(G,n)$ with a new top element adjoined, is the lift matroid of $G$ of size $n$ in the sense of Zaslavsky \cite{Bias2}. The full subsemigroup lattice of $B(G,n)$ and its BRSC are thus natural extensions of this important class of matroids and will be of interest to study by the methods of the current paper.

There are natural series of monoids parameterized by an integer $n$ and either a group $G$ or a field $F$. These include the full and partial transformation monoids, $T_n$ and $PT_n$ respectively, the monoid $M_n(F)$ of all $n \times n$ matrices over a field $F$ and the aforementioned Brandt semigroups $B(G,n)$. We ask if there are other natural series of semigroups and suggest to study their subsemigroup and congruence complexes.

\section*{Acknowledgments}

The first author acknowledges support from the Binational Science Foundation (BSF) of the United States and Israel, grant number 2012080. The second author acknowledges support from the Simons Foundation.
 The third author was partially supported by CMUP (UID/MAT/00144/2013), which is funded by FCT (Portugal) with national (MEC) and European structural funds (FEDER), under the partnership agreement PT2020.

\bigskip

{\sc Stuart Margolis, Department of Mathematics, Bar Ilan University,
  52900 Ramat Gan, Israel}

{\em E-mail address:} margolis@math.biu.ac.il

\bigskip

{\sc John Rhodes, Department of Mathematics, University of California,
  Berkeley, California 94720, U.S.A.}

{\em E-mail addresses}: rhodes@math.berkeley.edu, BlvdBastille@gmail.com

\bigskip

{\sc Pedro V. Silva, Centro de
Matem\'{a}tica, Faculdade de Ci\^{e}ncias, Universidade do
Porto, R. Campo Alegre 687, 4169-007 Porto, Portugal}

{\em E-mail address}: pvsilva@fc.up.pt

\end{document}